\documentclass{ article}
\usepackage{indentfirst}
\usepackage{amsfonts}

\usepackage{amssymb}
\usepackage{amsmath}
\usepackage{amsthm}

\usepackage{stmaryrd}
\usepackage{dsfont}
\usepackage{xfrac}
\usepackage{mathrsfs}

\begin{document}
\vspace*{.5cm}
\begin{center}
{\Large{\bf  Generalized almost para-contact manifolds}}\\
\vspace{.5cm}
 { Bayram \d{S}ahin and Fulya \d{S}ahin} \\
\end{center}

\vspace{.5cm}
\begin{center}
{ Department of Mathematics, Inonu University, 44280, Malatya-Turkey.}
\end{center}
\begin{center}
(emails:bayram.sahin@inonu.edu.tr, fulya.sahin@inonu.edu.tr)
\end{center}
\vspace{.5cm}

\noindent {\bf Abstract.} {\small  In this paper, we introduce generalized almost para-contact manifolds and obtain normality  conditions in terms of classical tensor fields. We show that such manifolds naturally carry certain Lie bialgebroid/quasi-Lie algebroid structures on them and we relate this new generalized manifolds with classical almost para-contact manifolds. The paper contains several examples.}\\

\section*{\bf 1.~Introduction}
\renewcommand{\theequation}{1.\arabic{equation}}
 As a unification and extension of usual notions of complex manifolds and symplectic manifolds, the notion of generalized complex manifolds was introduced by Hitchin \cite{Hitchin}. This subject has been studied widely in \cite{gualt} by Gualtieri. In particular, the notion of a generalized K\"{a}hler structure was introduced and studied by Gualtieri in the context of the theory of generalized geometric structures.  Later such manifolds and their submanifolds  have been studied in many papers. For instance, in \cite{Crainic}, Crainic gave necessary and sufficient conditions for a generalized almost complex manifold to be generalized complex manifold in terms of classical tensor fields. Generalized pseudo-K\"{a}hler structures have been also studied recently in \cite{Davidov}. Generalized geometry has received a reasonable amount of interest due to possible several relations with mathematical physics. Indeed, generalized K\"{a}hler
structures describe precisely the bi-Hermitian geometry arising in super-symmetric $\sigma-$ models \cite{GHR}.\\

   A central idea in generalized
geometry is that $TM\oplus TM^{\ast}$ should be thought of as a
generalized tangent bundle to manifold $M$. If $X$ and $\xi$ denote
a vector field and  a dual vector field on $M$ respectively, then we
write $(X,\alpha)$ (or $X+\alpha$) as a typical element of $TM\oplus
TM^{\ast}$. The space of sections of the vector bundle $TM\oplus TM^{\ast}$  is endowed with two natural $\mathbb{R}-$ bilineer operations: for the sections $(X,\alpha),(Y,\beta)$
of $TM\oplus TM^{\ast}=\mathcal{TM}$, a symmetric bilinear form $<,>$ is defined by
\begin{equation}
<X+\alpha,Y+\beta>=\frac{1}{2}(i_X\beta+i_Y\alpha),\label{eq:1.0}
\end{equation}
and the Courant bracket of two sections  is defined by
\begin{equation}
\llbracket(X,\alpha),(Y,\beta)\rrbracket=[X,Y]+L_{X}\beta-L_{Y}\alpha
\ -\frac{1}{2}d(i_{X}\beta-i_{Y}\alpha),\label{eq:1.1}
\end{equation}
where $d$, $L_{X}$ and $i_{X}$ denote exterior derivative, Lie
derivative and interior derivative with respect to $X$,
respectively. The Courant bracket is antisymmetric but, it does not
satisfy the Jacobi identity. In this paper we adapt the notions

\begin{equation}
\beta(\pi^{\sharp}\alpha)=\pi(\alpha,\beta)\quad\mathrm{and}\quad \omega_{\sharp}(X)(Y)=\omega(X,Y)\label{eq:1.01}
\end{equation}
 which are defined as
$\pi^{\sharp}:TM^{\ast} \rightarrow TM$,
$\omega_{\sharp}:TM\rightarrow TM^{\ast}$  for any 1-forms $\alpha$
and $\beta$, 2-form $\omega$ and bivector field $\pi$, and vector
fields $X$ and $Y$. Also we denote by $[,]_{\pi}$, the bracket on
the space of 1-forms on $M$ defined by
\begin{eqnarray}
[\alpha,\beta]_{\pi}=L_{\pi^{\sharp}\alpha}\beta-L_{\pi^{\sharp}\beta}\alpha-d\pi(\alpha,\beta).\label{eq:1.2}
\end{eqnarray}

 As an analogue of generalized complex structures on even dimensional
 manifolds, the concept of  generalized almost subtangent manifolds
 were introduced in \cite{vaisman} and such manifolds have been
 studied in \cite{vaisman} and  \cite{Wade1}. On the other hand, the notion of a generalized contact pair on a manifold $M$ was introduced by Poon and Wade in \cite{Poon-Wade}, see also \cite{Vaisman}, \cite{Vaisman2} and \cite{Wade2}. As we mention above, the framework of generalized almost complex structures puts almost
symplectic structures and almost complex structures on an equal footing. Similarly the notion of
a generalized almost contact structure unifies almost cosymplectic structures and almost contact
structures.\\

In this paper, we introduce  generalized almost  para-contact
structures/ manifolds and show that such manifolds include
para-contact manifolds as a subclass. Then we investigate normality
conditions for generalized almost para-contact manifolds in terms of
classical tensor fields and obtain certain Lie algebroid structures ( Courant
algebroid, quasi-Lie bialgebroid) on such manifolds. We give various
examples and show that classical almost para-contact manifolds can
be described in this generalized geometry.\\

The paper is organized as follows.  In section 2,
we recall some basic notions needed for the paper. In section 3, we
define generalized almost para-contact manifold, give examples and
obtain normality conditions for generalized almost para-contact
manifolds in terms of classical tensor fields. In section 4, we
investigate Lie algebroid structrures on generalized almost
para-contact manifolds. For this aim, we construct several
subbundles of big tangent bundle and show that all these subbundles
are isotropic and then we use this information to obtain a
characterization for quasi-Lie bialgebroid structure on a almost
generalized para-contact manifold. We also introduce to the notion
of strong generalized para-contact manifold and provide an example.
In section 5, we show that classical normal almost para-contact
manifolds are strong generalized para-contact manifolds. 

\section*{\bf 2.~Preliminaries}
\setcounter{equation}{0}
\renewcommand{\theequation}{2.\arabic{equation}}
An $(2n+1)$ dimensional smooth manifold $M$ has an almost para-contact structure $(\varphi,E,\eta)$ if it admits a tensor field $\varphi$ of type $(1,1)$, a vector field $E$ and a $1-$ form $\eta$ satisfying the following compatibility conditions

\begin{eqnarray}
\varphi(E)=0&,& \eta \circ \varphi=0\label{eq:-1.1}\\
\eta(E)=1&,&\varphi^2=id-\eta\otimes E\label{eq:-1.2}
\end{eqnarray}
and the tensor field $\varphi$ induces an almost paracomplex structure on each fibre on the distribution $\mathcal{D}$ generated by $\eta$. An immediate consequence of the definition of the almost para-contact structure is that the
endomorphism $\varphi$ has rank $2n$.\\

Let $M^{(2n+1)}$  be an almost para-contact manifold with structure $(\varphi, E, \eta)$ and consider the
manifold $M^{(2n+1)}\times \mathbb{R}$. We denote a vector field on $M^{(2n+1)}\times \mathbb{R}$ by $(X, f \frac{d}{dt} )$, where $X$ is
tangent to $M^{(2n+1)}$, $ t$ is the coordinate on $\mathbb{R}$, and $f$ is a $C^{\infty}$ function on $M^{(2n+1)}\times \mathbb{R}$. An
almost paracomplex structure $J$ on $M^{(2n+1)}\times \mathbb{R}$ is defined  by
$$J(X, f\frac{d}{dt})=(\varphi X + fE,\eta(X)\frac{d}{dt}).$$
If $J$ is integrable which means that the Nijenhuis tensor of $J$, $N_J$, vanishes, we say that the almost para-contact structure $(\varphi, E, \eta)$ is normal. For details, see:\cite{Kaneyuki-Willams} and \cite{Zamkovoy}.\\

We now recall the notion of generalized almost para-comlex structure on $\mathcal{TM}=TM\oplus TM^*$. A generalized almost paracomplex structure on $M$ is a vector bundle automorphism $\mathcal{J}$ on $\mathcal{TM}$ such that $\mathcal{J}^{2}= I$, $\mathcal{J}\neq I$ and $\mathcal{J}$ is orthogonal with respect to $<,>$, i.e.,
\begin{equation}
<\mathcal{J}e_1,e_2>+<e_1,\mathcal{J}e_2>=0,\, e_1,e_2\in \Gamma(\mathcal{TM}).\label{eq:2.0}
\end{equation}
A generalized almost paracomplex structure can be represented by
classical tensor fields as follows:
\begin{eqnarray}
\mathcal{J}= \left[ {\begin{array}{cc}
 a & \pi^{\sharp}  \\ \label{eq:2.1}
 \theta_{\sharp} & -a^\ast  \\
 \end{array} } \right]
\end{eqnarray}
where $\pi$ is a bivector on $M$, $\theta$ is a 2-form on $M$, $a :
TM \rightarrow TM$ is a bundle map, and $a^{\ast} : TM^{\ast}
\rightarrow TM^{\ast}$ is dual of $a$, for almost para-complex
structures see:\cite{vaisman} and \cite{Wade1}.\\

A  generalized almost  paracomplex structure is called integrable (or
just paracomplex structure) if $\mathcal{J}$ satisfies the following
condition
\begin{eqnarray}
 \llbracket \mathcal{J}\alpha,\mathcal{J}\beta\rrbracket-
\mathcal{J}(\llbracket \mathcal{J}\alpha, \beta\rrbracket + \llbracket \alpha,\mathcal{J}\beta\rrbracket)+ \llbracket \alpha,\beta\rrbracket
= 0, \label{eq:2.2}
\end{eqnarray}
 for all sections $\alpha,\beta \in \Gamma(\mathcal{TM})$.\\

\section*{\bf 3.~Generalized almost para-contact manifolds}
\setcounter{equation}{0}
\renewcommand{\theequation}{3.\arabic{equation}}
In this section, we are going to define generalized almost para-contact structure and obtain the normality conditions. We first propose the following definition.\\

\noindent{\bf Definition~1.~} A generalized almost para-contact structure $(\mathcal{F}, \mathcal{Z},\mathcal{\xi})$ on a smooth odd dimensional manifold $M$ consists of a bundle endomorphism $\mathcal{F}$ from $\mathcal{TM}$ to itself and a section $\mathcal{Z}+\mathcal{\xi}$ of $\mathcal{TM}$ such that
\begin{eqnarray}
\mathcal{F}+\mathcal{F}^*=0&,&\mathcal{\xi}(\mathcal{Z})=\mathcal{I}\label{eq:3.1}\\
\mathcal{F}(\mathcal{\xi})=0&,&\mathcal{F}(\mathcal{Z})=0\label{eq:3.2}\\
\mathcal{F}^2&=&\mathcal{I}-\mathcal{Z}\odot\mathcal{\xi}\label{eq:3.3}
\end{eqnarray}
where $\mathcal{Z}\odot\mathcal{\xi}$ is defined by
\begin{equation}
(\mathcal{Z}\odot\mathcal{\xi})(X+\alpha)=\mathcal{\xi}(X)\mathcal{Z}+\alpha(\mathcal{Z})\mathcal{\xi}.\label{eq:3.4}
\end{equation}
The bundle map $\mathcal{F}:\mathcal{TM}\longrightarrow \mathcal{TM}$ is given by
\begin{eqnarray}
\mathcal{F}= \left[ {\begin{array}{cc}
 F & \pi^{\sharp}  \\ \label{eq:3.5}
 \sigma_{\sharp} & -F^\ast  \\
 \end{array} } \right],
\end{eqnarray}
where $F:TM\longrightarrow TM$ is a bundle map. The following example shows that generalized almost para-contact structure is really a generalization of almost para-contact structure.\\

\noindent{\bf Example~1.~}Associated to any almost para-contact structure, we have a generalized almost para-contact structure
by setting
\begin{eqnarray}
\mathcal{F}= \left[ {\begin{array}{cc}
 \varphi & 0  \\ \nonumber
 0 & -\varphi^\ast  \\
 \end{array} } \right],
\end{eqnarray}
with the given vector field $E$ and 1-form $\eta$.\\

We give another example of generalized almost para-contact manifolds.\\

\noindent{\bf Example~2.~} Let $ H_3$ be the three-dimensional Heisenberg group and  $\{X_1,X_2,X_3\}$ a basis for its algebra
$\mathfrak{h}_3$ so that $[X_1,X_2] =- X_3$. Let $\{\alpha^1,\alpha^2, \alpha^3\}$ be a dual frame. Then $d\alpha^3 = \alpha^1\wedge \alpha^2$. Now for some real number $\vartheta$, we define
\begin{eqnarray}
F=\cosh\, \vartheta (X_2\otimes \alpha^2+X_3\otimes \alpha^3)&,&\mathcal{\xi}=\alpha^1,\,\mathcal{Z}=X_1\nonumber\\
\sigma=\sinh\,\vartheta (\alpha^2\wedge \alpha^3)&,&\pi=\sinh\, \vartheta (X_2\wedge X_3).\nonumber
\end{eqnarray}
We also , as given in (\ref{eq:3.5}), define
 \begin{eqnarray}
\mathcal{F}= \left[ {\begin{array}{cc}
 F & \pi^{\sharp}  \\ \nonumber
 \sigma_{\sharp} & -F^\ast  \\
 \end{array} } \right].
\end{eqnarray}
Then $(\mathcal{\xi},\mathcal{Z},\pi,\sigma,F)$ is a  generalized almost para-contact structure on $ H_3$.\\

We now obtain the normality conditions for a generalized almost para-contact structure. For this aim, we make the following definition.\\

\noindent{\bf Definition~2.~}A generalized almost para-complex structure $\mathcal{J}$ on $M\times \mathbb{R}$ is said
to be $M-$ adapted if it has the following three properties
\begin{enumerate}
  \item [(i)] $\mathcal{J}$ is invariant by
translation along $\mathbb{R}$.
  \item [(ii)] $\mathcal{J}(T\mathbb{R}\oplus 0)\subseteq 0\oplus TM^*$.
  \item [(iii)] $\mathcal{J}(0\oplus T\mathbb{R}^*)\subseteq TM\oplus 0$.
\end{enumerate}
 The invariance of $\mathcal{J} $ by translations means that the Lie derivatives $\frac{d}{dt}$  of
the classical tensor fields of $\mathcal{ J} $ defined by (\ref{eq:2.1}) vanish. If conditions (ii) and (iii)
are also imposed, it follows that the classical tensor fields of an $M-$ adapted,
generalized almost para-complex structure are of the form
\begin{equation}
a=F,\,  \pi=P+\mathcal{Z}\wedge \frac{d}{d t},\, \theta=\sigma+\mathcal{\xi}\wedge dt\label{eq:3.6}
\end{equation}
where $P$ is a bivector on $M$, $\sigma$ is a 2-form on $M$, $F :
TM \rightarrow TM$ is a bundle map. A generalized, almost para-contact structure will be called normal if the corresponding
$M-$ adapted, generalized almost para-complex structure on $M \times \mathbb{R}$ is integrable. The following theorem gives necessary and sufficient conditions for a generalized almost para-contact manifold.\\

\noindent{\bf Theorem~3.1.~} {\it A generalized, almost para-contact structure is normal if and only if the following conditions are satisfied.
\begin{enumerate}
  \item [(A1)] $P$ satisfies the equation
  \begin{equation}
  [P^{\sharp}\alpha_1,P^{\sharp}\beta_1]=P^{\sharp}([\alpha_1,\beta_1]_{P}).\label{eq:3.7}
  \end{equation}
  \item [(A2)] $P$ and $F$ are related by the following two formulas
  \begin{eqnarray}
  FP^{\sharp}&=&P^{\sharp}F^{*}\label{eq:3.8}\\
  F^*([\alpha_1,\beta_1]_P)&=&L_{P^{\sharp}\alpha_1}F^*\beta_1-L_{P^{\sharp}\beta_1}F^*\alpha_1+dP(\beta_1,F^*\alpha_1).\label{eq:3.9}
  \end{eqnarray}
  \item [(A3)] $P$, $\sigma$ and $F$ are related by the following four formulas
  \begin{eqnarray}
  i_Z\sigma=0, i_{\mathcal{\xi}}P=0&,&F^2=Id-P^{\sharp}\sigma^{\flat}-\mathcal{Z}^t\otimes\mathcal{\xi}\label{eq:3.10}\\
  N_F(X,Y)&=&P^{\sharp}(i_{X\wedge Y}(d\sigma^{\flat}))-(d\mathcal{\xi}(X,Y))\mathcal{Z}.\label{eq:3.11}
  \end{eqnarray}
  \item [(A4)] $\mathcal{\xi}$, $\mathcal{Z}$,$F$ and $\sigma^{\flat}$ are related the following formulas
  \begin{eqnarray}
  F(\mathcal{Z})=0, \mathcal{\xi} \circ F=0&,&(L_{FX}\mathcal{\xi})Y-(L_{FY}\mathcal{\xi})X=0\label{eq:3.12}\\
  d\sigma_F(X,Y,Z)&=&d\sigma(FX,Y,Z)+d\sigma(X,FY,Z)\nonumber\\
  &&+d\sigma(X,Y,FZ).\label{eq:3.13}
    \end{eqnarray}
  \item[(A5)] $\mathcal{\xi}$, $P$, $F$ and $\mathcal{Z}$ satisfy the following equations
  \begin{equation}
  L_{\mathcal{Z}}\mathcal{\xi}=0,\,L_{\mathcal{Z}}P=0,\, L_{\mathcal{Z}}F=0,\,L_{\mathcal{Z}}\sigma^{\flat}=0,\, L_{P^{\sharp}\alpha_1}\mathcal{\xi}=0\label{eq:3.14}
  \end{equation}
\end{enumerate}
for $X,Y,Z\in \Gamma(TM)$ and $\alpha_1, \beta_1 \in \Gamma(TM^*)$.}\\

\noindent{\bf Proof.~} The condition $\mathcal{J}^2=\mathcal{I}$ implies the first equations appearing in (A2), (A3) and (A4). For remaining parts, we need to check the integrability of $\mathcal{J}$, i.e.
$$\mathcal{N}(\mathcal{X},\mathcal{Y})=\llbracket \mathcal{J}\mathcal{X},\mathcal{J}\mathcal{Y}\rrbracket-\mathcal{J}\llbracket \mathcal{X},\mathcal{J}Y\rrbracket-\mathcal{J}\llbracket \mathcal{J}X,\mathcal{Y}\rrbracket+\llbracket \mathcal{X},\mathcal{Y}\rrbracket=0$$
for $\mathcal{X}, \mathcal{Y}, \mathcal{Z} \in \Gamma(\mathcal{TM}\oplus\mathbb{R})$. First, if $\mathcal{X}=(0,0, \alpha_1,0)$ and $\mathcal{Y}=(0,0,\beta_1,0)$ are $1-$ forms, then we have, as the vector field part,
$$[P^{\sharp}\alpha_1,P^{\sharp}\beta_1]=P^{\sharp}(L_{P^{\sharp}\alpha_1}\beta_1-L_{P^{\sharp}\beta_1}\alpha_1-dP(\alpha_1,\beta_1))$$
which is (\ref{eq:3.7}) due to (\ref{eq:1.2}).  For the $1-$ form part, we get (\ref{eq:3.9}).  Now for  vector fields $\mathcal{X}=(X, 0, 0,0)$ and $\mathcal{Y}=(Y,0,0,,0)$ we have
\begin{eqnarray}
[FX,FY]-F[FX,Y]-F[X,FY]+[X,Y]&-&P^{\sharp}(L_X\sigma^{\flat}(Y)-L_Y\sigma^{\flat}(X)\nonumber\\
+d\sigma(X,Y))-(L_Y\mathcal{\xi}(X))(\mathcal{Z})&+&(L_Y\mathcal{\xi}(X))(\mathcal{Z})=0.\nonumber
\end{eqnarray}
But using the following formula
$$i_{X\wedge Y}(d\sigma)=L_X(i_Y(\sigma))-L_Y(i_X(\sigma))+d(i_{X\wedge Y}\sigma)-i_{[X,Y]}\sigma$$
we arrive at
$$
N_F(X,Y)=P^{\sharp}(i_{X\wedge Y}(d\sigma^{\flat}))-(d\mathcal{\xi}(X,Y))\mathcal{Z}$$
which is (\ref{eq:3.11}). For $1-$ form part, we have
\begin{eqnarray}
L_FX\sigma^{\flat}(Y)-L_FY\sigma^{\flat}(X)&-&\frac{1}{2}d\sigma(Y,FX)+\frac{1}{2}d\sigma(X,FY)\nonumber\\
+F^*(L_X\sigma^{\flat}(Y)-L_Y\sigma^{\flat}(X)&+&d\sigma(X,Y))-\sigma^{\flat}([X,FY])\nonumber\\
&&-\sigma^{\flat}([X,FY])=0\label{eq:3.15})
\end{eqnarray}
and
\begin{equation}
FX(\mathcal{\xi}(Y))-FY(\mathcal{\xi}(X))-\mathcal{\xi}([FX,Y])-\mathcal{\xi}([X,FY])=0.\label{eq:3.16}
\end{equation}
Now using the following formula
\begin{eqnarray}
d\sigma(X,Y,Z)&=&L_X\sigma(Y,Z)+L_Y\sigma(Z,X)+L_Z\sigma(X,Y)\nonumber\\
&&-\sigma([X,Y],Z)-\sigma([Z,X],Y)-\sigma([Y,Z],X)\nonumber
\end{eqnarray}
then, (\ref{eq:3.15}) and (\ref{eq:3.16}) give (\ref{eq:3.12}) and (\ref{eq:3.13}). For $\mathcal{X}=(X,0,0,0)$ and $\mathcal{Y}=(0,0,\alpha_1,0)$, we have, as $1-$ form part,
\begin{eqnarray}
L_X\alpha_1-L_{FX}F^*\alpha_1-L_{P^{\sharp}\alpha_1}\sigma^{\flat}(X)&+&F^*(L_{FX}\alpha_1-L_XF^*\alpha_1)\nonumber\\
+\frac{1}{2}d(\alpha_1(F^2X)+\sigma(X,P\sharp^{\sharp}\alpha_1)&-&\alpha_1(X)+\alpha_1(\mathcal{Z})\mathcal{\xi}(X))\nonumber\\
&&-\sigma^{\flat}([X,P^{\sharp}\alpha_1])=0.\nonumber
\end{eqnarray}
Then from the third equation of (\ref{eq:3.10}) we arrive at
\begin{eqnarray}
L_X\alpha_1-L_{FX}F^*\alpha_1-L_{P^{\sharp}\alpha_1}\sigma^{\flat}(X)&+&F^*(L_{FX}\alpha_1-L_XF^*\alpha_1)\nonumber\\
-d\alpha_1(P^{\sharp}\sigma^{\flat}(X))-\sigma^{\flat}([X,P^{\sharp}\alpha_1])=0.\label{eq:3.17}
\end{eqnarray}
And for the vector field part, we obtain
$$
[FX,P^{\sharp}\alpha_1]-F[X,P^{\sharp}\alpha_1]-P^{\sharp}(L_{FX}\alpha_1-L_XF^*\alpha_1)=0
$$
which is equivalent to (\ref{eq:3.9}). And (\ref{eq:3.17}) is equivalent to (\ref{eq:3.15}). For $\mathcal{X}=(0,\frac{\partial }{\partial  t},0,0)$ and $\mathcal{Y}=(0,0,0,dt)$, we obtain the first equation in (\ref{eq:3.14}). In a similar way, we derive the other equations.\\

Now by using $F, \sigma, P$, we can write the endomorphism $\mathcal{F}$ of $\mathcal{TM}$  as
\begin{eqnarray}
\mathcal{F}\left( {\begin{array}{c}
 X  \\
 \alpha   \\
 \end{array} } \right)= \left[ {\begin{array}{cc}
 F & P^{\sharp}  \\ \label{eq:3.18}
 \sigma_{\sharp} & -F^\ast  \\
 \end{array} } \right]\left( {\begin{array}{c}
 X  \\ \label{eq:3.18}
 \alpha   \\
 \end{array} } \right).
\end{eqnarray}
Also using $\mathcal{Z}$ and $\mathcal{\xi}$, we can define the endomorphism $\mathfrak{Z}$ of $\mathcal{TM}$ given by
\begin{eqnarray}
\mathfrak{Z}\left( {\begin{array}{c}
 X  \\
 \alpha   \\
 \end{array} } \right)= \left[ {\begin{array}{cc}
 \mathcal{Z}\otimes\mathcal{\xi} & 0   \\ \label{eq:3.19}
 0 & (\mathcal{Z}\otimes\mathcal{\xi})^t  \\
 \end{array} } \right]\left( {\begin{array}{c}
 X  \\ \label{eq:3.18}
 \alpha   \\
 \end{array} } \right).
\end{eqnarray}
Now we give a characterization for normality in terms of $(\mathcal{F}, \mathfrak{Z}, (\mathcal{Z},\mathcal{\xi}))$.\\

\noindent{\bf Theorem~3.2.~}{\it Let $M$ be an odd dimensional manifold. Then $(\mathcal{F}, \mathfrak{Z}, (\mathcal{Z},\mathcal{\xi}))$ defines generalized normal contact structure on $M$ if and only if the following conditions are satisfied.}\\
\begin{enumerate}
  \item [(i)] {\it For $(X,\alpha)\in \Gamma(\mathcal{TM})$, $\mathcal{F}$ and $\mathfrak{Z}$ are related by the following formulas}
  \begin{equation}
  \mathcal{F}^2=\mathcal{I}-\mathfrak{Z}\label{eq:3.20}
  \end{equation}
  \begin{equation}
  \mathcal{N}_{\mathcal{F}}((X,\alpha),(Y,\beta))=\mathfrak{Z}\llbracket(X,\alpha),(Y,\beta)\rrbracket.\label{eq:3.21}
  \end{equation}
    \item [(ii)] {\it For $(X,\alpha) \in \Gamma(im \mathcal{F})$, $\mathcal{Z}$, $\mathcal{\xi}$ and $\mathcal{F}$ are related by the following formulas}
    \begin{equation}
  \parallel \mathcal{Z}\oplus \mathcal{\xi}\parallel_{g}=1\label{eq:3.22}
  \end{equation}
  \begin{equation}
  \llbracket(Z,0), \mathcal{F}((X,\alpha)\rrbracket=\mathcal{F}(L_{\mathcal{Z}}X,L_{\mathcal{Z}}\alpha),\label{eq:3.23}
  \end{equation}
  {\it where $g$ is the neutral metric of $\mathcal{TM}$ given
  by}
  $$g((X,\alpha),(Y,\beta))=\alpha(Y)+\beta(X).$$
    \item [(iii)]  {\it For $(X,\alpha) \in \Gamma(im \mathcal{F})$, $\mathcal{Z}$, $\mathcal{\xi}$ and $\mathcal{F}$ are related by the following formulas}
    \begin{equation}
  \flat_g \circ \mathcal{F}+\mathcal{F}^t \circ \flat_g=0,\mathcal{F} \circ \mathfrak{Z}=0\label{eq:3.24}
  \end{equation}
  \begin{equation}
  \llbracket \mathcal{F}(X,\alpha), (0, \mathcal{\xi})\rrbracket=\mathcal{F}(0,L_{X}\mathcal{\xi}).\label{eq:3.25}
  \end{equation}

\end{enumerate}

\noindent{\bf Proof.~} Using (\ref{eq:3.6}) and $\mathcal{J}^2=\mathcal{I}$ we have
\begin{eqnarray}
&&(F^2+P^{\sharp}\sigma^{\flat})(X)+(FP^{\sharp}-P^{\sharp}F^*)(\alpha)\nonumber\\
&&gF\mathcal{Z}+\mathcal{\xi}(X)\mathcal{Z}+(\sigma^{\flat}(X)-F^*\alpha+f\mathcal{\xi})(\mathcal{Z})\frac{d}{d t}\nonumber\\
&&(\sigma^{\flat} \circ F-F^*\circ \sigma^{\flat})(X)+(\sigma^{\flat}\circ P^{\sharp}+{F^*}^2)(\alpha)\nonumber\\
&&g\sigma^{\flat}(Z)-fF^*\mathcal{\xi}+\alpha(\mathcal{Z})\mathcal{\xi}+(\mathcal{\xi}(FX+P^{\sharp}\alpha+gZ))dt\nonumber\\
&&=X+f\frac{d}{d t}+\alpha+gdt.\label{eq:3.26}
\end{eqnarray}
Then (\ref{eq:3.20}), (\ref{eq:3.22}) and (\ref{eq:3.24}) follow from (\ref{eq:3.26}). On the other hand
we have
\begin{eqnarray}
\mathcal{J}(X+f\frac{d}{d t},\alpha+gdt)&=&\mathcal{F}(X,\alpha)+g(\mathcal{Z},0)+f(0,\mathcal{\xi})\nonumber\\
&&(\alpha(\mathcal{Z})\frac{d}{d t},\mathcal{\xi}(X)dt).\label{eq:3.27}
\end{eqnarray}
It is clear that $(X,\alpha) \in \Gamma(im\mathcal{F})$ if and only if $\mathcal{\xi}(X)=0$ and $\alpha(\mathcal{Z})=0$. Now we check the normality conditions. Since
$$\mathcal{TM}=im\mathcal{F}\oplus Sp\{(\mathcal{Z},0)\}\oplus Sp\{(0,\mathcal{\xi})\}\oplus Sp\{(\frac{d}{d t},0)\}\oplus Sp\{(0,dt)\},$$
it will be enough to consider combinations $(X,\alpha)\in \Gamma(im \mathcal{F})$, $(0,\mathcal{\xi})$, $(\frac{d}{d t},0)$, $(0,dt)$. For $(\frac{d}{d t},0)$, $(0,dt)$, since $\mathcal{J}(\frac{d}{d t},0)=(0,\mathcal{\xi})$ and $\mathcal{J}(0,dt)=(\mathcal{Z}, 0)$, we obtain
\begin{equation}
\mathcal{N}_{\mathcal{J}}((\frac{d}{d t},0), (0,dt))=0 \Leftrightarrow L_Z\mathcal{\xi}=0.\label{eq:3.28}
\end{equation}
For $(\mathcal{Z},0)$, $(0,\mathcal{\xi})$, we get
\begin{equation}
\mathcal{N}_{\mathcal{J}}((\mathcal{Z},0),(0,\mathcal{\xi}))=0 \Leftrightarrow L_Z\mathcal{\xi}=0.\label{eq:3.29}
\end{equation}
For $(X,\alpha)\in \Gamma(im\mathcal{F})$, $(0,dt)$, using (\ref{eq:3.27}) we have
\begin{equation}
\mathcal{N}_{\mathcal{J}}((X,\alpha), (0,dt))=0\Leftrightarrow \llbracket \mathcal{F}(X,\alpha),(\mathcal{Z},0)\rrbracket=\mathcal{F}(L_X\mathcal{Z},\mathcal{L}_{\mathcal{Z}}\alpha).\label{eq:3.30}
\end{equation}
For $(X,\alpha)\in \Gamma(im\mathcal{F})$,$(\frac{d}{d t},0)$, using (\ref{eq:3.27}) we obtain
\begin{equation}
\mathcal{N}_{\mathcal{J}}((X,\alpha), (\frac{d}{d t},0))=0 \Leftrightarrow \llbracket \mathcal{F}(X,\alpha),(0,\mathcal{\xi})\rrbracket=\mathcal{F}(0,L_X\mathcal{\xi}).\label{eq:3.31}
\end{equation}
For $(X,\alpha), (Y,\beta) \in \Gamma(im\mathcal{F})$, again from (\ref{eq:3.27}) we derive
\begin{equation}
\mathcal{N}_{\mathcal{J}}((X,\alpha), (Y,\beta))=0\Leftrightarrow \mathcal{N}_{\mathcal{F}}((X,\alpha),(Y,\beta))-\mathfrak{Z}\llbracket(X,\alpha),(Y,\beta)\rrbracket=0.\label{eq:3.32}
\end{equation}
Remaining combinations give  the conditions obtained in (\ref{eq:3.28})-(\ref{eq:3.31}). Note that the condition $L_Z\mathcal{\xi}=0$ in (\ref{eq:3.28}) is equivalent to (\ref{eq:3.21}).
\section*{\bf 4.~Lie bialgebroid structures on generalized almost para-contact manifolds}
\setcounter{equation}{0}
\renewcommand{\theequation}{4.\arabic{equation}}
In this section we are going to investigate necessary and sufficient
conditions for the existence of Lie bialgebroid structures on a generalized
para-contact manifold. We first recall some notions needed for this
section. A Lie algebroid structure on a real
 vector  bundle $A$ on a manifold $M$ is defined by a vector bundle
 map $\rho_{A}:A\rightarrow TM$, the anchor of $A$, and an
 $\mathbb{R}$-Lie algebra bracket on $\Gamma(A), [,]_{A}$ satisfying
 the Leibnitz rule
 \begin{center}
 $[\alpha,f \beta]_{A}=f[\alpha,\beta]_{A}+L_{\rho_{A}(\alpha)}(f)\beta$
 \end{center}
 for all $\alpha,\beta \in \Gamma(A), f \in C^{\infty}(M)$, where $L_{\rho_{A}(\alpha)}$ is
 the Lie derivative with respect to the vector field
 $\rho_{A}(\alpha)$. Suppose that $A\longrightarrow P$ is a Lie algebroid, and that its dual bundle
$A^*\longrightarrow P$ also carries a Lie algebroid structure. Then $(A;A^*)$ is a Lie bialgebroid
if for any $X, Y\in \Gamma(A)$,
$$d_*[X; Y ] = L_Xd_*Y- L_Y d_*X,$$
where $d_*$ is the  exterior derivative associated to $A^*$. On the other hand, a Courant algebroid $E\longrightarrow M$ is a vector bundle on $M$ equipped with a non-degenerate symmetric bilinear form $<,>$, a vector bundle map $\rho:E\longrightarrow TM$, a bilinear bracket $\circ$ on $\Gamma(E)$ satisfying
\begin{enumerate}
  \item [(c1)] $e_1\circ(e_2\circ e_3)=(e_1\circ e_2)\circ e_3+e_2 \circ (e_1\circ e_3)$,
  \item [(c2)] $e \circ e=\rho^*d<e,e>$,
  \item [(c3)] $L_{\rho(e)}<e_1,e_2>=<e \circ e_1,e_2>+<e_1,e\circ e_2>$,
  \item [(c4)] $\rho(e_1\circ e_2)=\rho[e_1,e_2]$,
  \item [(c5)] $e_1\circ fe_2=f(e_1\circ e_2)+L_{\rho(e_1)f}e_2$
\end{enumerate}
for all $e,e_1,e_2,e_3 \in \Gamma(E)$ and $f\in
C^{\infty}(M,\mathbb{R})$. It is known that $(TM\oplus TM^*)$ with
the non-degenerate metric $<,>$ given in (\ref{eq:1.0}) and Courant
bracket form a Courant algebroid. The anchor is the natural
projection  from the direct sum to the summand $ TM$. The following
result is well known.\\

\noindent{\bf Theorem~4.1.~}\cite{LWX}~{\it If $(A,A^*)$ is a Lie bialgebroid, then  $A\oplus A^*$ together with $([X,Y]_A,\rho_A, \llbracket,\rrbracket)$ is a Courant algebroid, where $\llbracket,\rrbracket$ denotes the Courant bracket given in (\ref{eq:1.1}). Conversely, in a Courant algebroid $(E,\rho,[,]_A,\circ)$, suppose that $L_1$ and $L_2$ are Dirac subbundles transversal to each other, then $(L_1,L_2)$ is a Lie bialgebroid.}\\

A real, maximal isotropic sub-bundle $L\subset TM\oplus TM^*$ is called an
almost Dirac structure. If $L$ is involutive, then the almost Dirac structure is said to be integrable,
or simply a Dirac structure. Similarly, a maximal isotropic and involutive para-complex sub-bundle
$L \subset (TM\oplus TM^*)_{\mathbb{C}}$ is called a para-complex Dirac structure.\\

A quasi-Lie bialgebroid \cite{Roytenberg} is a Lie algebroid $(A, [,]_A, \rho)$ equipped
with a degree-one derivation $d_*$ of the Gerstenhaber algebra $(\Gamma(\wedge^{\bullet}),[,]_A,\wedge)$ and
a $3-$ section of $A$, $X_A\in \Gamma(\wedge^3A)$ such that $d_*X_A=0$ and $d^2_*=[X_A,]$. If $X_A$ is the null section, then $d_*$ defines a structure of Lie algebroid on $A^*$ such that $d_*$ is a derivation of $[ , ]_A$. Let $(A, [,]_A, \rho)$ be a Lie algebroid and consider any closed 3-form
$\phi$.  Equipping $A^*$ with the null Lie algebroid structure, $(A^*,d_A, \pi)$ is canonically a
quasi-Lie bi-algebroid. We note that Lie algebroid structures have been also extended to Poisson-Nijenhuis structures \cite{GU}\\

Let $M$ be a generalized almost para-contact manifold and $\mathcal{F}$ the bundle map given in (\ref{eq:3.5}). Then it has one real eigenvalue, namely, $0$. The corresponding eigenbundle is trivialized by $\mathcal{Z}$ and $\mathcal{\xi}$, and we denote these bundles by $L_{\mathcal{Z}}$ and $L_{\mathcal{\xi}}$, respectively. Let $ker\mathcal{\xi}$ be the distribution on the manifold $M$ defined by the point-wise kernel of the 1-form
$\mathcal{\xi}$. Similarly, $ker \mathcal{Z}$ is the sub-bundle of $TM^*$ defined by the point-wise kernel of the vector field
$\mathcal{Z}$ with respect to its evaluation on differential 1-forms. On the paracomplexified bundle $(TM\oplus TM^*)_{\mathbb{C}}$, we have that $\mathcal{F}$ has three eigenvalues, namely, $0$, $1$ and $-1$. We now define
\begin{eqnarray}
E^{(1,0)}=\{\mathcal{X}+e\mathcal{F}\mathcal{X}\mid \mathcal{X}\in \Gamma(ker \mathcal{Z}\oplus ker\mathcal{\xi})\}\nonumber\\
E^{(0,1)}=\{\mathcal{X}-e\mathcal{F}\mathcal{X}\mid \mathcal{X}\in \Gamma(ker \mathcal{Z}\oplus ker\mathcal{\xi})\}\nonumber
\end{eqnarray}
where $e^2=1$, The extension $\mathcal{F}$ of the endomorphism $\mathcal{F}$ with eigenvalues
$\mp1$ to $(TM\oplus TM^*)_{\mathbb{C}}$ has eigenvalues $\mp e$, see:\cite{AMT} for para-complex structures on a real vector space and other notions. Then $L_{\mathcal{Z}}\oplus L_{\mathcal{\xi}}$ is the $0-$ eigenbundle, $E^{(1,0)}$ is the $1-$ eigenbundle and $E^{(0,1)}$ is $(-1)-$ eigenbundle. We have natural decomposition $(TM\oplus TM^*)_{\mathbb{C}}=L_{\mathcal{Z}}\oplus L_{\mathcal{\xi}}\oplus E^{(1,0)}\oplus E^{(0,1)}$. We now have the following four different paracomplex vector bundle which play different roles
\begin{eqnarray}
L=L_{\mathcal{Z}}\oplus E^{(1,0)}&,&\bar{L}=L_{\mathcal{Z}}\oplus E^{(0,1)}\nonumber\\
L^*=L_{\mathcal{\xi}}\oplus E^{(0,1)}&,&{\bar{L}}^*=L_{\mathcal{\xi}}\oplus E^{(1,0)}.\nonumber
\end{eqnarray}
Since $L_{\mathcal{Z}}$ is a real line bundle, its paraconjugation is itself. Therefore paracomplex conjugation sends $L$ to $\bar{L}$. All of these bundles are independent of the choice of representatives of a generalized almost contact structure.\\

\noindent{\bf Lemma~4.1.~}{\it The bundles $L$, $L^*$, $\bar{L}$, ${\bar{L}^*}$, $E^{(1,0)}$ and $E^{(0,1)}$ are isotropic with respect to $<,>$ given in (\ref{eq:1.0})}

\noindent{\bf Proof.~} For $X+\alpha$, using (\ref{eq:3.5}) we have
$$\mathcal{F}(X+\alpha)=FX+\pi^{\sharp}\alpha+\sigma_{\sharp}(X)-F^*\alpha.$$
Now if $X+\alpha\in \Gamma(ker\mathcal{\xi}\oplus ker\mathcal{Z})$, from (\ref{eq:1.0}) we get
$$
<\mathcal{Z},\mathcal{F}(X+\alpha)>=\frac{1}{2}i_{\mathcal{Z}}\sigma_{\sharp}(X)-\alpha(F\mathcal{Z}).$$
Then using (\ref{eq:3.2}) and (\ref{eq:3.10}) we have
\begin{equation}
<\mathcal{Z},\mathcal{F}(X+\alpha)>=0. \label{eq:4.1}
\end{equation}
In a similar way we also have
\begin{equation}
<\mathcal{\xi},\mathcal{F}(X+\alpha)>=0. \label{eq:4.}
\end{equation}
Thus if $X+\alpha\in \Gamma(ker\mathcal{\xi}\oplus ker\mathcal{Z})$, then $\mathcal{F}(X+\alpha)\in \Gamma(ker\mathcal{\xi}\oplus ker\mathcal{Z})$. For $X+\alpha, Y+\beta \in \Gamma(ker\mathcal{\xi}\oplus ker\mathcal{Z})$, using (\ref{eq:3.5}) we have
\begin{eqnarray}
&&<X+\alpha+e\mathcal{F}(X+\alpha),Y+\beta+e\mathcal{F}(X+\beta)>=\beta(X)+e\beta(FX)\nonumber\\
&&+e\beta(\pi^{\sharp}\alpha)+e\sigma_{\sharp}(Y)(X)+e\sigma_{\sharp}(Y)(e(FX+\pi^{\sharp}\alpha))-eF^*\beta(X)\nonumber\\
&&-F^*\beta(FX+\pi^{\sharp}\alpha)+\alpha(Y)+e\alpha(FY+\pi^{\sharp}\beta)+e\sigma_{\sharp}(X)(Y)\nonumber\\
&&+e\sigma_{\sharp}(X)(e(FY+\pi^{\sharp}\beta))-eF^*\alpha(Y)-F^*\alpha(FY+\pi^{\sharp}\beta).\label{eq:4.3}
\end{eqnarray}
Since $\sigma$ is symmetric and
$$\sigma(Y,FX)+\sigma(X,FY)=F^*(\sigma_\sharp(Y)(X)+\sigma_\sharp(X)(Y))$$
we have
\begin{equation}
\sigma(Y,FX)+\sigma(X,FY)=0.\label{eq:4.4}
\end{equation}
On the other hand, using (\ref{eq:3.10}) and (\ref{eq:1.01}) we derive
\begin{equation}
-\alpha(F^2Y)+\sigma_{\sharp}(Y)(\pi^{\sharp}\alpha)+\alpha(Y)=0.\label{eq:4.5}
\end{equation}
Now putting (\ref{eq:4.4}) and (\ref{eq:4.5}) in (\ref{eq:4.3}) we obtain
$$<X+\alpha+e\mathcal{F}(X+\alpha),Y+\beta+e\mathcal{F}(X+\beta)>=0$$
which shows that $E^{(1,0)}$ is isotropic. Taking the paracomplex conjugation in above computation, we see that $E^{(0,1})$ is also isotropic. Also from (\ref{eq:4.1}) and (\ref{eq:4.}) we obtain that the pairings between $L_\mathcal{\xi}$ or $L_\mathcal{Z}$ with those $E^{(1,0})$ and $E^{(0,1})$ are equal to zero. Thus $L$ is isotropic. Similarly, we find that $L^*$ is isotropic.\\

We now present the following notion.\\

\noindent{\bf Definition~3.~}Given a generalized almost para-contact
structure, if the space $\Gamma(L)$ of sections of the associated
bundle $L$ is closed under the Courant bracket, then the generalized
almost para-contact structure is simply called a generalized
para-contact structure.\\

Thus by considering the notions of Dirac structures, quasi-Lie
bialgebroids and Definition~3, we have the following result.\\

\noindent{\bf Corollary~4.1.~}{\it When $\mathcal{J} =(\mathcal{Z},
\mathcal{\xi}, \pi, \sigma, F)$ represents a generalized
para-contact structure, the associated bundle $L$ is a Dirac
structure. In addition, the bundle $L^*$ is a transversal isotropic
complement of $L$ in the Courant algebroid $((TM \oplus
TM^*)_{\mathbb{C}},<,>, \llbracket,\rrbracket)$, where $\llbracket, \rrbracket$ denotes the Courant bracket. As a
result,   the pair $L$ and $L^*$ is a quasi-Lie bialgebroid.}\\

We now investigate Lie bialgebroid structures on a generalized
para-contact manifold.\\

\noindent{\bf Theorem~4.2.~}{\it Let $\mathcal{J} =(\mathcal{Z},
\mathcal{\xi}, \pi, \sigma, F)$ be  an integrable generalized para-contact structure. The pair $L$ and $L^*$ forms a
Lie bialgebroid if and only if $d\mathcal{\xi}$ is of type $(1, 1)$
with respect to the map $\mathcal{F}$ on $\Gamma((ker
\mathcal{Z}\oplus ker \mathcal{\xi})_C)$.}\\

\noindent{\bf Proof.~} We first note that  the inclusion of $L$ in
$(TM \oplus TM^*)_{\mathbb{C}}$ followed by the natural projection onto the
first summand is an anchor map. When $L$ is closed under the Courant
bracket, the restriction of the Courant bracket to $L$ completes the
construction of a Lie algebroid structure on $L$. Since $\mathcal{J}
=(\mathcal{Z}, \mathcal{\xi}, \pi, \sigma, F)$ is a generalized
para-contact structure, $L$ is closed under the Courant bracket.
From (\cite{gualt}, Proposition~3.27), we know that a maximal
isotropic sub-bundle $\tilde{E}$  of $(TM \oplus TM^*)_{\mathbb{C}}$ is involutive if
and only the Nijenhuis operator $Nij$ given by
\begin{equation}
Nij(A,B,C)=\frac{1}{3}\{<\llbracket A,B\rrbracket,C>+<\llbracket B,C\rrbracket,A>+<\llbracket C,A\rrbracket,B> \label{eq:4.6}
\end{equation}
vanishes on $\tilde{E}$. Since $L$ is closed with respect to the Courant
bracket, by conjugation we have
$$\llbracket \Gamma(E^{(0,1})),\Gamma(E^{(0,1}))\rrbracket \subseteq \Gamma(L_{\mathcal{Z}}\oplus
E^{(0,1}))=\Gamma(\bar{L}).$$  Isotropic $\bar{L}$ implies that
$<L_F,L_F\oplus E^{(0,1})>=<E^{(0,1}),L_F\oplus E^{(0,1})>=0$. Thus
if $\mathcal{X}_1$, $\mathcal{X}_2$ and
$\mathcal{X}_3$ are all sections of $E^{(0,1})$, then
$Nij(\mathcal{X}_1,\mathcal{X}_2,\mathcal{X}_3)=0$. Thus $Nij=0$ if and only if
$$Nij(\mathcal{X}_1,\mathcal{X}_2,\mathcal{\xi})=0$$
for $\mathcal{X}_1$, $\mathcal{X}_2 \in \Gamma(E^{(0,1)})$. Now for $X_1\in \Gamma(ker\mathcal{\xi})$ and $\mathcal{\alpha}_1\in \Gamma(ker\mathcal{Z})$, by taking
$\mathcal{X}_1=X_1+\alpha_1-e\mathcal{F}(X_1+\alpha_1)$ and $\mathcal{X}_2=X_2+\alpha_2-e\mathcal{F}(X_2+\alpha_2)$ we find
$$\mathcal{X}_1=X_1+\alpha_1-eFX_1-e\pi^{\sharp}\alpha_1-e\sigma_{\sharp}(X_1)+eF^{*}(\alpha_1)$$
and
$$
\rho(\mathcal{X}_1)=X_1-eFX_1-e\pi^{\sharp}\alpha_1,
$$
where $\rho:L^*\longrightarrow TM$ is the anchor map. Thus we have
$$Nij(\mathcal{X}_1,\mathcal{X}_2,\mathcal{\xi})=\frac{1}{3}\{<\rho\llbracket\mathcal{X}_1,\mathcal{X}_2\rrbracket,\mathcal{\xi}>+<\rho\llbracket \mathcal{X}_1,\mathcal{\xi}\rrbracket,\rho\mathcal{X}_2>+<\rho\llbracket\mathcal{\xi},\mathcal{X}_1\rrbracket,\rho\mathcal{X}_2>\}.$$
Hence we get
\begin{eqnarray}
Nij(\mathcal{X}_1,\mathcal{X}_2,\mathcal{\xi})&=&\frac{1}{6}\mathcal{\xi}([\rho\mathcal{X}_1,\rho\mathcal{X}_2])+\frac{1}{6}(L_{\rho(\mathcal{X}_2)}\mathcal{\xi})(\rho(\mathcal{X}_1))-\frac{1}{12}(d(i_{\rho(\mathcal{X}_2)}\mathcal{\xi})(\rho(\mathcal{X}_1))\nonumber\\
&&-\frac{1}{6}(L_{\rho(\mathcal{X}_1)}\mathcal{\xi})(\rho(\mathcal{X}_1))+\frac{1}{12}(d(i_{\rho(\mathcal{X}_1)}\mathcal{\xi})(\rho(\mathcal{X}_1)).\label{eq:4.7}
\end{eqnarray}
Since $X_2\in \Gamma(ker\mathcal{\xi})$, by using the second equation of (\ref{eq:3.10}) we derive
\begin{equation}
d(i_{\rho(\mathcal{X}_2)}\mathcal{\xi})(\rho(\mathcal{X}_1))=0.\label{eq:4.8}
\end{equation}
Using this in (\ref{eq:4.7}) we obtain
$$
Nij(\mathcal{X}_1,\mathcal{X}_2,\mathcal{\xi})=\frac{1}{6}\{\mathcal{\xi}([\rho\mathcal{X}_1,\rho\mathcal{X}_2])+(L_{\rho(\mathcal{X}_2)}\mathcal{\xi})(\rho(\mathcal{X}_1))-(L_{\rho(\mathcal{X}_1)}\mathcal{\xi})(\rho(\mathcal{X}_1))\}.
$$
Then exterior derivative and  (\ref{eq:4.8})  imply that
$$Nij(\mathcal{X}_1,\mathcal{X}_2,\mathcal{\xi})=\frac{1}{3}\{-d\mathcal{\xi}(\rho(\mathcal{X}_1),\rho(\mathcal{X}_2))+d\mathcal{\xi}(\rho(\mathcal{X}_2),\rho(\mathcal{X}_1))-d\mathcal{\xi}(\rho(\mathcal{X}_1),\rho(\mathcal{X}_2))\}.$$
Thus we arrive at
$$
Nij(\mathcal{X}_1,\mathcal{X}_2,\mathcal{\xi})=-d\mathcal{\xi}(\rho(\mathcal{X}_1),\rho(\mathcal{X}_2))
$$
which gives
\begin{equation}
Nij=-\mathcal{Z}\wedge \rho^*(d\mathcal{\xi})^{(2,0)},\label{eq:4.10}
\end{equation}
where $\rho^*(d\mathcal{\xi})^{(2,0)}$ is the $\wedge E^{(1,0)}-$ component of the pullback of the $d\mathcal{\xi}$ via the anchor map $\rho:L^*\longrightarrow TM$. Then proof follows from Theorem ~4.1, Lemma~4.1 and (\ref{eq:4.10}).\\

The above result motivates us to give the following definition.\\

\noindent{\bf Definition~4.~} An generalized almost para-contact structure is called a strong generalized
para-contact structure if both $L$ and $L^*$ are closed under the Courant bracket.\\

Here is an example for strong generalized para-contact structure.\\

\noindent{\bf Example~3.~} Let $ H_3$ be the three-dimensional Heisenberg group and  $\{X_1,X_2,X_3\}$ a basis for its algebra
$\mathfrak{h}_3$ so that $[X_1,X_2] =- X_3$. Let $\{\alpha^1,\alpha^2, \alpha^3\}$ be a dual frame. Then $d\alpha^3 = \alpha^1\wedge \alpha^2$. Now for $t=r\cosh\, \vartheta+er\sinh\,\vartheta$, we define
\begin{eqnarray}
F_t=\frac{2r\sinh\, \vartheta}{(1-r^2)} (X_2\otimes \alpha^2+X_3\otimes \alpha^3)&,&\mathcal{\xi}=\alpha^1,\,\mathcal{Z}=X_1\nonumber\\
\sigma_t=\frac{-r^2+2r\cosh\,\vartheta-1}{(1-r^2)}(\alpha^2\wedge \alpha^3)&,&\pi_t= \frac{r^2+2r\cosh\,\vartheta+1}{(1-r^2)}(X_2\wedge X_3).\nonumber
\end{eqnarray}
Then $(\mathcal{\xi},\mathcal{Z},\pi_t,\sigma_t,F_t)$ is a family of
generalized almost para-contact structures. There are two
subfamilies of this family, determined by $|t|^2 = r^2 < 1$ and
$|t|^2 = r^2 > 1$.
 The corresponding
bundle $L_t$ and $L^*_t$ are trivialized as follows;
\begin{eqnarray}
L_t&=&span\{X_1,X_2+e\alpha^3+e\mathcal{F}(X_2+e\alpha^3),X_3-e\alpha^2+e\mathcal{F}(X_3-e\alpha^2)\}\nonumber\\
&=&span\{X_1,-(r^2+r\cosh v,)X_2+er\sinh vX_2+e(-r^2+r\cosh v)\alpha^3\nonumber\\
&&-r\sinh v\alpha^3,-(r^2+r\cosh\, v)X_3+r\sinh\, v\alpha^2 +er\sinh\, vX_3\nonumber\\
&&+e(r^2-r\cosh\, v)\alpha^2\}\nonumber
\end{eqnarray}
\begin{eqnarray}
L^*_t&=&span\{\alpha^1,\alpha^2+eX_3-e\mathcal{F}(\alpha^2+eX_3),\alpha^3-eX_2-e\mathcal{F}(\alpha^3-eX_2)\}\nonumber\\
&=&span\{\alpha^1, -r\sinh\,vX_3+(-r^2+r\cosh\,v)\alpha^2-e(r^2+r\cosh\,v)X_3\nonumber\\
&&+er\sinh\,v\alpha^2, r\sinh\, vX_3+(-r^2+r\cosh\, v)\alpha^3+e(r^2+r\cosh\,v)X_3\nonumber\\
&&+r\sinh\,v\alpha^3\}.\nonumber
\end{eqnarray}
It is easy to see that all Courant brackets on $L^*_t$ are trivial. The only non-zero Courant bracket on $L_t$ is
\begin{eqnarray}
\llbracket X_1,-(r^2&+&r\cosh v\,)X_2+er\sinh vX_2+e(-r^2+r\cosh v)\alpha^3-r\sinh v\alpha^3\rrbracket\nonumber\\
&=&(r^2+r\cosh \,v-er\sinh\,v)X_3-(r\sinh\,v+e(r^2-r\cosh\,v))\alpha^2,\nonumber
\end{eqnarray}
which is belong to $L_t$. Thus $(\mathcal{\xi},\mathcal{Z},\pi_t,\sigma_t,F_t)$ is a family of strong generalized para-contact structures.

\section*{\bf 5.~Almost para-contact structures as an example of strong generalized para-contact structures}
\setcounter{equation}{0}
\renewcommand{\theequation}{5.\arabic{equation}}
In this section we are going to show that an almost para-contact
structure has a natural strong generalized para-contact structure on
it. Let $(\varphi,E,\eta)$ be an almost para-contact structure on a
manifold $M$. An almost para-contact structure is a normal almost
para-contact structure \cite{Zamkovoy} if
\begin{equation}
N_{\varphi}(X,Y)=2d\eta(X,Y)E, L_E\eta=0,
L_E\varphi=0,\label{eq:5.1}
\end{equation}
where $N_{\varphi}$ is defined by
\begin{equation}
N_{\varphi}(X,Y)=[\varphi X,\varphi Y]+\varphi^2[X,Y]-\varphi[\varphi X,Y]-\varphi[X,\varphi Y]\label{eq:5.2}
\end{equation}
for any vector fields $X, Y$. \\

\noindent{\bf Lemma~5.1.~}{\it  Let $(\varphi,E,\eta)$ be a normal almost para-contact structure on a manifold $M$. If $\mathcal{X}$ is a section of $E^{(1,0})$, then $\llbracket F,\mathcal{X} \rrbracket$ is again a section of $E^{(1,0)}$}.\\

 \noindent{\bf Proof.~}For $\mathcal{X}=X+\alpha \in \Gamma(ker
E\oplus ker \eta)$, using (\ref{eq:1.1}) we have
\begin{equation}
\llbracket E,\mathcal{X}+e\mathcal{F}\mathcal{X}
\rrbracket=L_EX+L_E\alpha+e(L_E\varphi
X-L_E\varphi^*\alpha)\label{eq:5.3}
\end{equation}
due to $X\in \Gamma(ker\eta)$ and $\alpha \in \Gamma(Ker E)$. Since
$M$ is a normal almost para-contact manifold we have (\ref{eq:5.1}).
Replacing $X$ by $E$ in (\ref{eq:5.1}) we have
$$\varphi^2[E,Y]=-\eta([E,Y])E.$$
Then we get
\begin{equation}
L_EY=\varphi L_E\varphi Y.\label{eq:5.4}
\end{equation}
(\ref{eq:5.1}) and (\ref{eq:5.4}) imply that
\begin{equation}
\varphi L_EY=L_E\varphi Y.\label{eq:5.5}
\end{equation}
Thus using (\ref{eq:5.5}) in (\ref{eq:5.3}) and considering dual almost para-contact structure $\varphi^*$, we get
$$\llbracket E,\mathcal{X}+e\mathcal{F}\mathcal{X} \rrbracket=L_EX+L_E\alpha+e(\varphi L_E X-\varphi^*L_E\alpha)$$
which proves assertion.\\

In the sequel we show that a normal almost para-contact manifold is actually a strong generalized para-contact manifold.\\

\noindent{\bf Theorem~5.1.~}{\it If $\mathcal{J} =(\mathcal{Z},
\mathcal{\xi}, \pi, \sigma, F)$ represents a generalized almost para-contact structure associated to
a classical normal almost para-contact structure on an odd-dimensional manifold $M$, then it is a
strong generalized para-contact structure.}\\

\noindent{\bf Proof.~} For $X, Y\in \Gamma(ker \eta)$, applying $\eta$ to (\ref{eq:5.2}) we have
\begin{equation}
\eta([\varphi X, \varphi Y])=2d\eta(X,Y). \label{eq:5.6}
\end{equation}
Since $\varphi X$ and $\varphi Y$ are also sections of $ker \eta$, we derive
\begin{equation}
-d\eta(X,Y)=d\eta(\varphi X, \varphi Y).\label{eq:5.7}
\end{equation}
Now for $X, Y\in \Gamma(ker \eta)$, we have
$$
\llbracket X+e\varphi X,Y+e\varphi Y \rrbracket=[X,Y]+[\varphi X, \varphi
Y]+e([\varphi X, Y]+[X, \varphi Y]).$$ Hence we get
\begin{eqnarray}
\llbracket X+e\varphi X,Y+e\varphi Y \rrbracket&=&[X,Y]+\varphi^2[\varphi X, \varphi Y]+\eta([\varphi X, \varphi Y])E\nonumber\\
&&+e([\varphi X, Y]+[X, \varphi Y]).\label{eq:5.8}
\end{eqnarray}
Since $M$ is a normal almost para-contact manifold, from (\ref{eq:5.1}), (\ref{eq:5.2}) and (\ref{eq:5.6}) we obtain
\begin{equation}
\eta([\varphi X, \varphi Y])E+[X,Y]+\varphi^2[\varphi X, \varphi
Y]=\varphi([\varphi X, Y]+[X, \varphi Y]).\label{eq:5.9}
\end{equation}
Putting (\ref{eq:5.9}) in (\ref{eq:5.8}) we arrive at
$$
\llbracket X+e\varphi X,Y+e\varphi Y\rrbracket=\varphi([\varphi X, Y]+[X, \varphi Y])+e([\varphi X, Y]+[X, \varphi Y]).$$
Hence we have
\begin{eqnarray}
\llbracket X+e\varphi X,Y+e\varphi YE \rrbracket&=&\varphi([\varphi X, Y]+[X, \varphi Y])\nonumber\\
&&+e(\varphi^2([\varphi X, Y]+[X, \varphi Y])+\eta([\varphi X,
Y]+[X, \varphi Y]E).\nonumber
\end{eqnarray}
On the other hand, from (\ref{eq:5.7}) we find
$$\eta([\varphi X, Y]+[X, \varphi Y])=0.$$
Using this in above equation, we get
\begin{eqnarray}
\llbracket X+e\varphi X,Y+e\varphi Y \rrbracket&=&\varphi([\varphi X, Y]+[X, \varphi Y])\nonumber\\
&&+e\varphi(\varphi([\varphi X, Y]+[X, \varphi Y]))\label{eq:5.10}
\end{eqnarray}
which shows that $\llbracket X+e\varphi X,Y+e\varphi Y \rrbracket$ is a
section of $E^{(1,0})$.  Now for $X\in \Gamma(ker \eta)$ and $\alpha
\in \Gamma(ker E)$, from (\ref{eq:1.1}), we obtain
\begin{equation}
\llbracket X+e\varphi X,\alpha-e\varphi^* \alpha
\rrbracket=L_X\alpha-L_{\varphi
X}\varphi^*\alpha-e(L_X\varphi^*\alpha-L_{\varphi
X}\alpha).\label{eq:5.11}
\end{equation}
Evaluating the  above expression on the Reeb field $E$ and using
(\ref{eq:5.1})-(\ref{eq:5.3}) we see that both parts of the above
expression are  sections of $ker E$. On the other hand, for any
vector field $Y \in \Gamma(ker \eta)$, we have
$$
-\varphi^*(L_X\varphi^*\alpha-L_{\varphi X}\alpha)(Y)=-(L_X\varphi^*)(\varphi Y)+(L_{\varphi X}\alpha)(\varphi Y).$$
Then observing that $\varphi^2X =X$ on $ker \eta$, we obtain

\begin{eqnarray}
-\varphi^*(L_X\varphi^*\alpha-L_{\varphi X}\alpha)(Y)&=&-X\alpha(Y)+\varphi X\alpha(\varphi Y)\nonumber\\
&&+\alpha(\varphi[X, \varphi Y]-[\varphi X, \varphi Y]).\nonumber
\end{eqnarray}
From (\ref{eq:5.1}) and (\ref{eq:5.2}) we have
$$[\varphi X, \varphi Y]-\varphi[X, \varphi Y]=\varphi[\varphi X,  Y]-[X,Y].$$
Using this expression in above equation, we derive
\begin{eqnarray}
-\varphi^*(L_X\varphi^*\alpha-L_{\varphi X}\alpha)(Y)&=&-X\alpha(Y)+\varphi X\alpha(\varphi Y)\nonumber\\
&&+\alpha([X, Y]-\varphi[\varphi X, Y]).\nonumber\\
&=&-(L_X\alpha)(Y)+(L_{\varphi X}\varphi^*\alpha)(Y).\nonumber
\end{eqnarray}
Hence we arrive at
\begin{equation}
-\varphi^*(L_X\varphi^*\alpha-L_{\varphi X}\alpha)=-L_X\alpha+L_{\varphi X}\varphi^*\alpha\label{eq:5.12}
\end{equation}
on $ker E$. Applying $\varphi^*$ to (\ref{eq:5.12}) we have
$$
-(L_X\varphi^*\alpha-L_{\varphi
X}\alpha)+(L_X\varphi^*\alpha-L_{\varphi
X}\alpha)(E)\eta=\varphi^*(-L_X\alpha+L_{\varphi X}\varphi^*\alpha)
$$
due to the transpose of the  formula in (\ref{eq:-1.2}). Since $M$
is a normal almost para-contact manifold, we obtain
$(L_X\varphi^*\alpha-L_{\varphi X}\alpha)(E)=0$, thus we have
\begin{equation}
(L_X\varphi^*\alpha-L_{\varphi
X}\alpha)=\varphi^*(L_X\alpha-L_{\varphi
X}\varphi^*\alpha).\label{eq:5.13}
\end{equation}
Using (\ref{eq:5.13}) in (\ref{eq:5.11}) we find
$$
\llbracket X+e\varphi X,\alpha-e\varphi^* \alpha
\rrbracket=L_X\alpha-L_{\varphi
X}\varphi^*\alpha-e\varphi^*(L_X\alpha-L_{\varphi X}\varphi^*\alpha)
$$
which shows that $\llbracket X+e\varphi X,\alpha-e\varphi^* \alpha
\rrbracket \in E^{(1,0})$. Since the Courant bracket between two forms
are zero, this shows that we have proved that the Courant bracket
between two sections of $E^{(1,0)}$ is again a section of
$E^{(1,0)}$. Since the bundle $E^{(1,0})$ is $E$ invariant, the
bundle is closed under the Courant bracket. Also we see that
(\ref{eq:5.7}) shows that $\rho^*(d\eta)^{(2,0)}=0$. Therefore $L^*$
is also closed with respect to the Courant bracket . Thus proof is
completed by Lemma 5.1 and Definition~4.


\begin{thebibliography}{99}
\bibitem{AMT}  D. V. Alekseevsky,  C. Medori and  A. Tomassini, Homogeneous para-K\"{a}hler Einstein manifolds, Russian Math. Surveys 64:1 (2009), 1-43.



\bibitem{Gualtieri-Cavalcanti}G. R. Cavalcanti and M. Gualtieri, Generalized complex geometry and T-duality, arxiv:1106.1747.

\bibitem {Crainic} M. Crainic, Generalized complex structures and Lie brackets, Bull. Braz. Math. Soc., New Series 42(4),
(2011),  559-578.

\bibitem{Davidov}J. Davidov, G. Grantcharov, O. Mushkarov,   M. Yotov, Generalized pseudo-K\"{a}hler structures, Comm. Math. Phys. 304(1), (2011), 49-68.

\bibitem{GHR}S. J. Gates Jr, C. M. Hull, M. R$\breve{o}$cek, Twisted multiplets and new supersymmetric nonlinear $\sigma-$ models,
Nuclear Phys. B,  248(1), (1984), 157-186.

\bibitem{GU} J. Grabowski and  P. Urbanski, Lie algebroids and Poisson Nijenhuis structures, Rep. Math. Phy. 40, (1997), 195-208.




\bibitem {gualt}M. Gualtieri, Generalized complex geometry, Ph.D. thesis, Univ. Oxford, arXiv:ma.th.DG/0401221, (2003).

\bibitem {Hitchin}N. Hitchin, Generalized Calabi-Yau manifolds, Q. J.
Math., 54, (2003), 281-308.





\bibitem{Kaneyuki-Willams}S. Kaneyuki,  F. L. Willams, Almost paracontact and parahodge structures on manifolds. Nagoya Math.
J. 99, (1985), 173-187.


\bibitem{LWX}Z-J. Liu, A. Weinstein, P. Xu,  Manin triples for Lie bialgebroids, J. Differential Geom. Volume 45(3), (1997), 547-574.










\bibitem{Poon-Wade}Y. S. Poon, A.  Wade, Generalized contact structures, J. Lond. Math. Soc.  83(2), (2011), 333-352.
\bibitem{Roytenberg}D. Roytenberg, Quasi-Lie bialgebroids and twisted Poisson manifolds. Lett. Math. Phys. 61
(2002), 123-137.




\bibitem{Vaisman}I. Vaisman, Dirac structures and generalized complex structures on $TM\times \mathbb{R}^h$  . Adv. Geom. 7(3), (2007), 453-474.


\bibitem {vaisman}I. Vaisman, Reduction and submanifolds of generalized complex manifolds, Dif. Geom. and its appl.,
25, (2007), 147-166.

\bibitem{Vaisman2}I. Vaisman,  From generalized K\"{a}hler to generalized Sasakian structures. J. Geom. Symmetry Phys. 18 (2010), 63-86.

\bibitem {Wade1}A. Wade, Dirac structures and paracomplex manifolds, C. R. Acad. Sci.
Paris, Ser. I, 338, (2004), 889-894.
\bibitem{Wade2}A. Wade,  Local structure of generalized contact manifolds. Differential Geom. Appl. 30(1), (2012), 124-135








\bibitem{Zamkovoy}S. Zamkovoy, Canonical connections on paracontact manifolds, Ann. Glob. Anal. Geom.,  36, (2009), 37-60.





\end{thebibliography}
\end{document}